\documentclass{article}    % Enable this line and disable the
\usepackage{graphicx}          % Include this line if your
\usepackage{epsfig}
\usepackage{graphicx}
\usepackage{makeidx}
\usepackage{multicol}
\usepackage{verbatim}
\usepackage{subfigure}
\usepackage{mathrsfs}
\usepackage{snapshot}

\usepackage{amsmath}

\usepackage{crop}
\crop
\usepackage{amssymb}%,minilisp,mycltxdefs}%,amsthm}%,verbatim}
\usepackage{graphicx}

\def\E{\mathbb{E}}

\newcommand{\zb}{\left(}
\newcommand{\zj}{\right)}

\newcommand{\sino}{\varepsilon^{{\rm (s)}}_n(\theta)}
\newcommand{\ino}{\varepsilon_n(\theta)}
\newcommand{\sdino}{\varepsilon^{{\rm (s)}}_{\theta n}(\theta)}

\newcommand{\sdinoT}{\varepsilon^{{\rm (s)T}}_{\theta n}(\theta)}

\newcommand{\sigs}{\sigma^2_n(\theta)}

\def\ba{\begin{array}}
\def\ea{\end{array}}
\def\bi{\begin{itemize}}
\def\ei{\end{itemize}}

\newtheorem{definition}{Definition}[section]

\newtheorem{theorem}{Theorem}[section]
\newtheorem{lemma}{Lemma}[section]

\begin{document}

%----
%\begin{frontmatter}
%\runtitle{Insert a suggested running title}  % Running title for regular
                                              % papers but only if the title
                                              % is over 5 words. Running title
                                              % is not shown in output.
\begin{comment}
----
\title{Identification of Finite Dimensional Linear Systems Driven by L\'evy processes} % Title, preferably not more
                                                % than 10 words.

%\thanks[footnoteinfo]{This paper was not presented at any IFAC
%meeting. Corresponding author M.~T.~Cicero. Tel. +XXXIX-VI-mmmxxi.
%Fax +XXXIX-VI-mmmxxv.}

\author[GL]{L\'aszl\'o Gerencs\'er}\ead{gerencser.laszlo@sztaki.mta.hu},    % Add the
\author[MM]{M\'at\'e M\'anfay}\ead{manfay@sztaki.mta.hu}            % e-mail address

\address[GL]{MTA SZTAKI}  % Please supply
\address[MM]{MTA SZTAKI, Central European University}             % full addresses
%\address[Baiae]{The White House, Baiae}        % here.

\begin{keyword}                           % Five to ten keywords,
linear systems, stochastic systems, L\'evy processes, system
identification, financial modelling              % chosen from the IFAC
\end{keyword}                             % keyword list or with the
                                          % help of the Automatica
                                          % keyword wizard
----
\end{comment}
\title{ECF identification of GARCH systems driven by L\'evy processes}
\date{}
\author{M\'at\'e M\'anfay and L\'aszl\'o Gerencs\'er and Zsanett Orlovits}

\maketitle
%----
%\begin{abstract}                          % Abstract of not more than 200 words.
L\'evy processes are widely used in financial mathematics,
telecommunication, economics, queueing theory and natural sciences
for modelling. We propose an essentially asymptotically efficient estimation method for the system parameters of general autoregressive conditional heteroscedasticity (GARCH) processes.
%Price processes are then defined as a corresponding geometric
%L\'evy process, implying the fact that returns are
%independent. In this paper we propose an alternative class of
%models allowing to describe dependence between return data.
As an alternative to
the maximum likelihood (ML) method we develop and analyze a novel
identification method by adapting the so-called empirical
characteristic function method (ECF) originally devised for
estimating parameters of c.f.-s from i.i.d. samples. Precise
characterization of the errors of these estimators will be given,
and their asymptotic covariance matrices will be obtained. 
%----
%\end{abstract}

%----
%\end{frontmatter}

\section{Basic properties of GARCH processes}

An important stylized fact of financial time series is that the conditional variance of the return process is not constant in time. This feature of financial data can be expressed by saying that it has a time varying volatility. Therefore the process cannot be modeled by linear systems. Thus, in particular, to analyze the dynamics of highly volatile financial instruments such as indices, foreign exchange rates and commodities, a more sophisticated model should be proposed that can reflect the dynamic volatility of past data resulting in the well-known phenomena of volatility clustering.

The first model that captured the above mentioned stylized fact was introduced by Engle \cite{arch_engle}. His model, the so called autoregressive conditional heteroscedasticity (ARCH) model, was  refined by Bollerslev \cite{garch_bollerslev}. Bollerslev's GARCH (generalized ARCH) model is one of the most widely accepted models recently in the area of financial modeling.

In this paper we tackle one of key problems of the statistical analysis of GARCH models, the parameter estimation problem, see \cite{gerencser2012real}. When it comes to the identification, the most principled method in the literature is the quasi-maximum likelihood method, see for example \cite{garch_berkes}. The main objective of this paper is to study the possibility of adapting the ECF method to GARCH processes with i.i.d. driving noise having known characteristic function. This possibility has not been attracted many researchers until recently. In \cite{garch_specification} a goodness of fit test is applied using the empirical characteristic function, while in \cite{garch_fourier} a Fourier type method is presented for power GARCH processes. Xu in \cite{garch_xu_normal} proposes to estimate the parameters of a GARCH model with normal driving noise using the ECF method and presents some empirical investigations.

Technically, the special type of a GARCH$(r,s)$ model to be studied in this paper is defined via the equations
\begin{align}
y_n&=\sigma_n \Delta L_n \label{eq:sigma*_def'}\\
\sigma_n^2-\gamma^*&=\sum_{i=1}^r \alpha_i^*(y^2_{n-i}-\gamma^*)+\sum_{j=1}^{s}\beta_j^*(\sigma^2_{n-j}-\gamma^*), \label{eq:sigma*_def}
\end{align}
where $ -\infty < n < + \infty.$ The driving noise $\Delta L_n$ is obtained as the
increment of a L\'evy process $(L_t)$ with $
-\infty < t < + \infty$, and $L_0=0$,
over an interval $[(n-1)h,nh),$ with $h>0$ being a fixed sampling
interval, and $ -\infty < n < + \infty$. The noise characteristic will be denoted by $\eta^*,$ i.e. the characteristic function of $\Delta L_n$ is $\varphi(u,\eta^*).$
We assume that $L_t$ has zero mean and ${\rm Var}(\Delta L_n)=1.$ Let $\mathscr{F}^{\Delta L}$ denote the natural filtration, i.e. $\mathscr{F}_{n}^{\Delta L}=\sigma \left\{ \Delta L_k: k \leq n \right\}$. Under the above conditions $\gamma^*$ is the conditional variance of $y_n$ and $\sigma_n$ given $\{y_i:i<n\}.$ The unknown parameter vector $\theta^*$ is defined as $
\theta^*=(\alpha^*_0,\alpha^*_1,\ldots,\alpha^*_r,\beta^*_1,\ldots,\beta^*_r)^T.
$
The second order properties of a GARCH process was given by Bollerslev, see \cite{garch_bollerslev}.
\begin{theorem}
The GARCH$(r,s)$ process defined by $(\ref{eq:sigma*_def})$ and $(\ref{eq:sigma*_def'})$ is second-order stationary with
$$
\E\left[ y_n\right]=0,~ {\rm Cov}(y_n,y_m)=0~for~n\neq m
$$
and
$$
\E\left[ y^2_n\right]=\E\left[ \sigma^2_n\right]=\frac{\alpha_0^*}{1-\sum_{i=1}^r \alpha_i^*-\sum_{j=1}^{s}\beta_j^*}
$$
if and only if
$$
\sum_{i=1}^r \alpha_i^*+\sum_{j=1}^{s}\beta_j^*<1.
$$
\end{theorem}
\begin{definition}
We say that a L\'evy process $(L_t)$ satisfies the moment condition of order $Q$ if
$$\int_{\mathbb{R}}|x|^q \nu(dx)<\infty$$
holds for $1 \leq q \leq Q,$ where the L\'evy measure of $(L_t)$ is denoted by $\nu(dx).$
\end{definition}
Define the polynomials
\begin{equation}
C^*(q^{-1})=\sum_{i=1}^r \alpha_i^* q^{-1} \text{ ~and~ } D^*(q^{-1})=1-\sum_{j=1}^s \beta_j^* q^{-1},
\end{equation}
with $q^{-1}$ being the backshift operator. In order to guarantee the invertibility of the sensitivity matrix we assume that $C^*$ and $D^*$ are relative prime. Using these polynomials $(\ref{eq:sigma*_def})$ can be written in the following compact form:
\begin{equation}
D^*(q^{-1})(\sigma_n^2-\gamma^*)=C^*(q^{-1})(y_n^2-\gamma^*).
\end{equation}
Let us define the $(r+s)$-dimensional state vector
\begin{equation}
\label{eq:state_garch}
X^*_n=(y_n^2,\ldots,y^2_{n-r+1},\sigma_n^2,\ldots,\sigma_{n-s+1}^2)^T.
\end{equation}
It is easy to check that the dynamics of $(X^*_n)$ is then
\begin{equation}
X^*_{n+1}=A_{n+1}^*X^*_n+u^*_{n+1}, \quad n \in \mathbb{Z},
\end{equation}
where $A_n^* \in \mathbb{R}^{(r+s)\times(r+s)}$ is defined in terms of $(\Delta L_n)$ as
%\begin{equation}$$
%$$\begin{pmatrix}
%  \alpha_1^* (\Delta L_n)^2 & \alpha_2^* (\Delta L_n)^2 & \cdots & \alpha_{r-1}^* (\Delta L_n)^2 & \alpha_r^* (\Delta L_n)^2 & \beta_1^*    (\Delta L_n)^2 & \beta_2^* (\Delta L_n)^2 & \ldots & \beta_{s-1}^* (\Delta L_n)^2 & \beta_s^* (\Delta L_n)^2  \\
%  1 & 0 & \cdots & 0 & 0 & 0 & 0 & \cdots & 0 & 0 \\
%  0 & 1 & \cdots & 0 & 0 & 0 & 0 & \cdots & 0 & 0 \\
%  \vdots & \vdots & \ddots & \vdots & \vdots & \vdots & \vdots & \ddots & \vdots & \vdots \\
%  0 & 0 & \cdots & 1 & 0 & 0 & 0 & \ldots & 0 & 0 \\
%  \alpha_1^*  & \alpha_2^*  & \cdots & \alpha_{r-1}^*  & \alpha_r^*  & \beta_1^*  & \beta_2^*  & \ldots & \beta_{s-1}^*  & \beta_s^*   \\
%  0 & 0 & \cdots & 0 & 0 & 1 & 0 & \cdots & 0 & 0 \\
%  0 & 0 & \cdots & 0 & 0 & 0 & 1 & \cdots & 0 & 0 \\
%  \vdots & \vdots & \ddots & \vdots & \vdots & \vdots & \vdots & \ddots & \vdots & \vdots \\
%  0 & 0 & \cdots & 0 & 0 & 0 & 0 & \cdots & 1 & 0 \\
% \end{pmatrix}$$
%\end{equation}
$$
A_n^*=\left(
        \begin{array}{cc}
          A_{n;1,1}^* & A_{n;1,2}^* \\
          A_{n;2,1}^* & A_{n;2,2}^* \\
        \end{array}
      \right),
$$
where
\begin{align}
A_{n;1,1}^*&=\left(
              \begin{array}{ccccc}
                \alpha_1^* (\Delta L_n)^2 & \alpha_2^* (\Delta L_n)^2 & \cdots & \alpha_{r-1}^* (\Delta L_n)^2 & \alpha_r^* (\Delta L_n)^2 \\
                1 & 0 & \cdots & 0 & 0 \\
                0 & 1 & \cdots & 0 & 0 \\
                \vdots & \vdots & \ddots & \vdots & \vdots \\
                0 & 0 & \cdots & 1 & 0 \\
              \end{array}
            \right) \\
A_{n;1,2}^*&=\left(
              \begin{array}{ccccc}
                \beta_1^*    (\Delta L_n)^2 & \beta_2^* (\Delta L_n)^2 & \ldots & \beta_{s-1}^* (\Delta L_n)^2 & \beta_s^* (\Delta L_n)^2  \\
                0 & 0 & \cdots & 0 & 0 \\
                0 & 0 & \cdots & 0 & 0 \\
                \vdots & \vdots & \ddots & \vdots & \vdots \\
                0 & 0 & \ldots & 0 & 0 \\
              \end{array}
            \right)
            \end{align}
            \begin{align}
A_{n;2,1}^*&=\left(
              \begin{array}{ccccc}
                \alpha_1^*  & \alpha_2^*  & \cdots & \alpha_{r-1}^*  & \alpha_r^* \\
                0 & 0 & \cdots & 0 & 0 \\
                0 & 0 & \cdots & 0 & 0 \\
                \vdots & \vdots & \ddots & \vdots & \vdots \\
                0 & 0 & \ldots & 0 & 0 \\
              \end{array}
            \right) \\
A_{n;2,2}^*&=\left(
              \begin{array}{ccccc}
                \beta_1^*  & \beta_2^*  & \ldots & \beta_{s-1}^*  & \beta_s^* \\
                1 & 0 & \cdots & 0 & 0 \\
                0 & 1 & \cdots & 0 & 0 \\
                \vdots & \vdots & \ddots & \vdots & \vdots \\
                0 & 0 & \cdots & 1 & 0 \\
              \end{array}
            \right)
\end{align}
and
$$
u_n^*=(\alpha_0^* (\Delta L_n)^2,0,\ldots,0,\alpha_0^*,0,\ldots,0)^T,
$$
for each $n,$ with $\alpha_0^*=\gamma^*\left(1-\sum_{i=1}^r \alpha^*_i-\sum_{j=1}^s \beta^*_j\right)$.
Note that $(A_n^*,u_n^*), n\in \mathbb{Z}$ is a sequence of i.i.d. random matrices. Moreover, $(X_n^*)$ is a Markov process with unobservable components. The above given state space representation, which is the slight modification of the one introduced by Bougerol and Picard \cite{Bougerol_Picard}, will be useful for proving $L$-mixing properties of $(y_n),(\sigma_n)$ and related processes.

Define the $p=r+s+1$-dimensional parameter vector
$$
\theta=(\alpha_0,\alpha_1,\ldots,\alpha_r,\beta_1,\ldots,\beta_r)^T
$$
and the real domain
$$
D=\left\{ \theta \left| \sum_{i=1}^r \alpha_i+\sum_{j=1}^s \beta_j <1\right.\right\}.
$$
Note that for $\theta \in D$ the corresponding GARCH process with parameter vector $\theta$ is well-defined. We also define a corresponding complex domain $D_{\epsilon}$ by
$$
D_{\epsilon}=\left\{ \theta \left| \operatorname{Re}\left(\sum_{i=1}^r \alpha_i+\sum_{j=1}^s \beta_j\right) <1,~|\operatorname{Im}(\alpha_i)|<\epsilon,~|\operatorname{Im}(\beta_j)|<\epsilon  \right.  \right\}
$$

 Let $D_{\epsilon}^* \subset {\rm int~ } D_{\epsilon}$ be a compact domain such that $\theta^* \in {\rm int~ } D_{\epsilon}^*.$ For a fixed value $\theta \in D_{\epsilon}$ we invert the GARCH system to recover the driving noise. Define the process $(\sigma_n(\theta))$ in terms of $y_n$:
\begin{equation}\label{eq:sigma_inverse}
\sigma_n^2 (\theta)-\gamma=\sum_{i=1}^r \alpha_i(y^2_{n-i}-\gamma)+\sum_{j=1}^{s}\beta_j(\sigma^2_{n-j}(\theta)-\gamma),
\end{equation}
with initial values $y_n=0,~\sigma_n^2(\theta)=\gamma, \text{ for all } n \leq 0.$
Then the estimated driving noise is defined as
\begin{equation}\label{eq:vareps_def}
\varepsilon_n(\theta)=\frac{y_n}{\sigma_n(\theta)}
\end{equation}
for $n \geq 0.$ Note that for $\theta=\theta^*$ the stationary solution of the inverse is
$$ \varepsilon^{(\rm s)}_n(\theta^*)=\frac{y_n}{\sigma^{{\rm (s)}}_n(\theta)}=\Delta L_n,$$
which is obtained by letting
$-\infty<n<\infty.$ Note that if $\theta=\theta^*,$ then $\sigma_n(\theta)$ recovers $\sigma_n$ at least in a statistical sense.

\section{ML method for GARCH processes}\label{sec:garch_ML}

In this section we develop and implement a ML estimator for GARCH models driven by L\'evy processes. As in the subsequent sections we will present an ECF identification method corresponding to the third stage of three-stage identification method for linear systems that is essentially asymptotically efficient. Before proceeding to this in this section we compute the asymptotic covariance of the ML method itself, to be able to compare the two asymptotic covariance matrices.

It can be shown along the lines of the ML method for linear systems that for the joint density function $f_Y$ of $\left( y_1,\ldots,y_n\right)$ and the joint density function $f_{\Delta L}$ of $\left( \Delta L_1,\ldots,\Delta L_n\right)$ we have
$$
f_Y \left( y_1,\ldots,y_n\right)=\prod_{k=1}^n \sigma_k(\theta)^{-1} f_{\Delta L} \left( \Delta L_1,\ldots,\Delta L_n \right),
$$
 because the determinant of the Jacobian of the transformation $\left( \Delta L_1,\ldots,\Delta L_n\right)\rightarrow \left( y_1,\ldots,y_n\right)$ is $\prod_{k=1}^n \sigma_k(\theta)^{-1}.$ We have that
the asymptotic cost function of the ML estimator is given by
\begin{equation}
\begin{split}
W(\theta)=\lim_{n\rightarrow \infty}\E\left[ - \log \left( f\left( \varepsilon_n(\theta)\right) \sigma_n(\theta)^{-1}\right) \right]=
\E\left[ - \log \left( f\left( \varepsilon^{{\rm (s)}}_n(\theta)\right) \sigma^{{\rm (s)}}_n(\theta)^{-1}\right) \right]= \\
\E\left[ - \log \left( f\left( \varepsilon^{{\rm (s)}}_n(\theta)\right)\right) +\log \left( \sigma^{{\rm (s)}}_n(\theta)\right)\right],
\end{split}
\end{equation}
where $f$ denotes the density function of $\Delta L_n.$
It is easy to check that $W_{\theta}(\theta^*)=0$ holds. For,
\begin{equation}
\begin{split}
W_{\theta}(\theta^*)= \E \left[ \frac{f'(\Delta L_n)}{f(\Delta L_n)} \Delta L_n \frac{\sigma^{{\rm (s)}}_{\theta n}(\theta^*)}{\sigma^{{\rm (s)}}_{n}(\theta^*)}+\frac{\sigma^{{\rm (s)}}_{\theta n}(\theta^*)}{\sigma^{{\rm (s)}}_{n}(\theta^*)}\right]= \\
\E\left[\E \left[ \left.\left(\frac{f'(\Delta L_n)}{f(\Delta L_n)} \Delta L_n +1\right)\frac{\sigma^{{\rm (s)}}_{\theta n}(\theta^*)}{\sigma^{{\rm (s)}}_{n}(\theta^*)}\right| \mathscr{F}^{\Delta L}_{n-1}\right]\right]=
\E \left[ \frac{\sigma^{{\rm (s)}}_{\theta n}(\theta^*)}{\sigma^{{\rm (s)}}_{n}(\theta^*)} \E \left[ \frac{f'(\Delta L_n)}{f(\Delta L_n)} \Delta L_n+1\right]\right],
\end{split}
\end{equation}
because $\varepsilon_{\theta n}(\theta^*)=-\Delta L_n \frac{\sigma_{\theta n}(\theta^*)}{\sigma_n(\theta^*)},$ and
$\sigma^{{\rm (s)}}_{\theta n}(\theta^*)$ and $\sigma^{{\rm (s)}}_{n}(\theta^*)$ are $\mathscr{F}^{\Delta L}_{n-1}$ measurable. Note that under appropriate regularity conditions on $f$ we have
\begin{equation}
\label{eq:E=-1}
\E \left[ \frac{f'(\Delta L_n)}{f(\Delta L_n)} \Delta L_n\right]=\int_{\mathbb{R}} f'(x) x dx=[xf(x)]_{-\infty}^{\infty}-\int_{\mathbb{R}} f(x) dx=-1,
\end{equation}
which implies the claim.
For the Hessian of $W$ write
\begin{equation}
\begin{split}
&W_{\theta \theta}(\theta)=-\lim _{n \rightarrow \infty}\left(\E \left[ \frac{f''(\sino)}{f(\sino)}\sdino \sdinoT -\frac{f'^{2}(\sino)}{f^2(\sino)}\sdino \sdinoT+\frac{f'(\sino)}{f(\sino)}\varepsilon^{{\rm (s)}}_{\theta \theta n}(\theta) \right]\right. \\
&+\left.\E\left[ \frac{\sigma^{{\rm (s)}}_{\theta \theta n}(\theta)\sigma^{{\rm (s)}}_n(\theta)-\sigma^{{\rm (s)}}_{\theta n}(\theta)\sigma^{{\rm (s)T}}_{\theta n}(\theta)}{\sigma^{{\rm (s)}2}_n(\theta)}\right]\right).
\end{split}
\end{equation}
Using again the fact that $\sigma^{{\rm (s)}}_{\theta n}(\theta^*)$ and $\sigma^{{\rm (s)}}_{n}(\theta^*)$ are $\mathscr{F}^{\Delta L}_{n-1}$ measurable we get that $W_{\theta \theta} (\theta^*)$ equals to
\begin{equation}
\begin{split}
&\E \left[ \frac{\sigma^{{\rm (s)}}_{\theta n}(\theta^*)\sigma^{{\rm (s)T}}_{\theta n}(\theta^*)}{\sigma^{{\rm (s)}2}_n(\theta^*)}\right] \E \left[ \frac{f'^2(\Delta L_n)}{f^2(\Delta L_n)}(\Delta L_n)^2- \frac{f''(\Delta L_n)}{f(\Delta L_n)}(\Delta L_n)^2-2\frac{f'(\Delta L_n)}{f(\Delta L_n)}\Delta L_n-1 \right]+ \\
&\E \left[ \frac{\sigma^{{\rm (s)}}_{\theta \theta n}(\theta^*)}{\sigma^{{\rm (s)}}_n(\theta^*)}\right] \E\left[ \frac{f'(\Delta L_n)}{f(\Delta L_n)}\Delta L_n+1\right].
\end{split}
\end{equation}
Using that under appropriate technical conditions we have
$$ \int_{\mathbb{R}}f''(x)x^2 dx=[f'(x)x^2]_{-\infty}^{\infty}-2\int_{\mathbb{R}}f'(x)x dx=2,$$ the previous formula can be written as
\begin{equation}
\E \left[ \frac{\sigma^{{\rm (s)}}_{\theta n}(\theta^*)\sigma^{{\rm (s)T}}_{\theta n}(\theta^*)}{\sigma^{{\rm (s)}2}_n(\theta^*)}\right] \E \left[ \frac{f'^2(\Delta L_n)}{f^2(\Delta L_n)}(\Delta L_n)^2-1\right].
\end{equation}
By almost identical calculation we get that the covariance of the gradient of log-likehood function
$$l_{\theta,n}(\theta)=\frac{f'(\Delta L_n)}{f(\Delta L_n)} \varepsilon^{{\rm (s)}}_{\theta n}(\theta)+\frac{\sigma^{{\rm (s)}}_{\theta n}(\theta)}{\sigma^{{\rm (s)}}_{n}(\theta)}$$
at $\theta=\theta^*$ is given by
\begin{equation}
{\rm Cov}(l_{\theta,n}(\theta^*),l^T_{\theta,n}(\theta^*))=\E \left[ \frac{\sigma^{{\rm (s)}}_{\theta n}(\theta^*)\sigma^{{\rm (s)T}}_{\theta n}(\theta^*)}{\sigma^{{\rm (s)}2}_n(\theta^*)}\right] \E \left[ \frac{f'^2(\Delta L_n)}{f^2(\Delta L_n)}(\Delta L_n)^2-1\right].
\end{equation}
Thus we have the following lemma.
\begin{lemma} Let $\hat{\theta}_N$ the ML estimate of the parameters of a GARCH process.
Then the asymptotic covariance of $\sqrt{N}\left( \hat{\theta}_N-\theta^*\right)$ is
\begin{equation}
\mu^{-1} (M^*)^{-1},
\end{equation}
with
$$ M^*=\E\left[\frac{\sigma^{{\rm (s)}}_{\theta N}(\theta^*)\sigma^{{\rm (s)T}}_{\theta N}(\theta^*)}{\sigma^{{\rm (s)}2}_N(\theta^*)}\right],$$
and
\begin{equation}\label{eq:mu_def2}
\mu=\E \left[ \frac{f'(\Delta L_N)^2}{f^2(\Delta L_N)}(\Delta L_N)^2-1\right].
\end{equation}
\end{lemma}
%The precise definition of the random vector $\hat{\theta}_N$ is analogous with that of the solution of $(\ref{eq:eta_felder}).$
A very nice interpretation of this $\mu$ is that it can be also obtained as the Fisher information of a scale parameter estimation problem. Suppose that we are given an i.i.d. realization of the scaled random variable $\lambda \Delta L_1,$ with the true value of $\lambda$ being $\lambda^*=1.$ Then the $\lambda$-dependent density $f(x,\lambda)$ of $\lambda \Delta L_1$ is
$$
f(x,\lambda)=f\left(\frac{x}{\lambda}\right)\frac{1}{\lambda},
$$
where $f(\cdotp)$ denotes the density function of $\Delta L_1.$
Write
\begin{equation}
\frac{\partial}{\partial \lambda} \log f(x,\lambda)=\frac{f'\zb \frac{\lambda}{x}\zj}{f\zb\frac{\lambda}{x}\zj}\zb-\frac{x}{\lambda^2}\zj-\frac{1}{\lambda},
\end{equation}
Hence taking into account (\ref{eq:E=-1}) the Fisher information reads as
\begin{equation}
\begin{split}
\E\left[ \left.\left(\frac{\partial}{\partial \lambda} \log f(\Delta L_1,\lambda)\right)^2\right|_{\lambda=1}\right]=\E\left[ \zb -\frac{f'(\Delta L_1)}{f(\Delta L_1)}\Delta L_1-1 \zj^2\right]= \\
\E \left[ \frac{f'^{2}(\Delta L_1)}{f^2(\Delta L_1)}\zb \Delta L_1 \zj^2+1+2\frac{f'(\Delta L_1)}{f(\Delta L_1)}\Delta L_1\right]=\E \left[ \frac{f'^{2}(\Delta L_1)}{f^2(\Delta L_1)}\zb \Delta L_1 \zj^2-1\right].
\end{split}
\end{equation}
%$$
%\E\left[ \frac{f'(\Delta L_1)}{f(\Delta L_1)}\Delta L_1\right]=\int_{\mathbb{R}}f'(x)x dx=\left[ f(x)x\right]_{-\infty}^{\infty}-\int_{\mathbb{R}}f(x)dx=-1.
%$$
Therefore we get the following lemma.
 \begin{lemma} $\mu$ in (\ref{eq:mu_def2}) can be interpreted as a Fisher information of a scale parameter estimate.
 \end{lemma}
 In analogy with the analysis of the efficiency of the three-stage method for linear systems this property of $\mu$ will have a key role in proving the essentially asymptotic efficiency of the ECF method for GARCH systems.

\section{ECF method for GARCH processes}\label{sec:ECFmethod_GARCH}

Now we turn to the problem of identifying the parameters of a GARCH process by adapting the approach of the ECF method. The ideas presented in this section show several similarities with those of \cite{automatica_own}, yet we will see that the different model structure poses numerous new problems. Despite the fact that the dynamics of a GARCH process can be described as a Markov process, the method presented in \cite{markov_ecf} does not solve this problem as it is not capable of dealing with unobservable components. For GARCH models only $(y_n)$ is observable and $(\sigma_n)$ is a latent process. The paper of Carrasco, Chernov, Florens and Ghysels \cite{markov_ecf} tackles the problem of estimating the parameters of an observable Markov process. Hereby we briefly summarize their findings. Let $X_t$ be a Markov process that is generated with some unknown parameter vector $\theta_0.$ Let $\varphi(s | X_t; \theta)$ denote the conditional characteristic function
$$
\E \left[ e^{isX_{t+1}} | X_t \right].
$$
The score functions used in the method are defined by
$$
h(r,s,X_t,X_{t+1};\theta)=e^{irX_t}\left( e^{isx_{t+1}} -\varphi\left(s |X_t;\theta_0 \right)\right).
$$
They prove that under some conditions using continuum moment condition yields an estimator that reaches the Cramer-Rao bound.

While this is a very attractive result, it does not solves the problems we consider in this thesis. The process $X_t$ is supposed to be observable, their proposed method cannot handle latent components. The presence of latent component is natural in GARCH processes, hence the method is not applicable for such processes.
For such non-Markovian processes they propose to use the joint characteristic function instead of the conditional one. Define the joint characteristic function as
$$
\varphi(\tau,t,n)=\E \left[ e^{i\tau^T Y_t} \right],
$$
where $\tau=(\tau_0,\ldots,\tau_n)^T$ and $Y_t=(X_t,\ldots,X_{t+n}).$
This problem has been considered in \cite{feuerverger_JCF} by Feuerverger. The unknown parameter $\theta$ is estimated by fitting the empirical joint characteristic function to the joint characteristic function using a weighting function. Feuerverger showed that this estimator is as efficient as the one that obtains $\hat{\theta}$ by solving
$$
\frac1N \sum_{k=1}^{N} \frac{\partial}{\partial \theta} \ln f(X_{k+n} | X_{k+n-1},\ldots,X_{k+n};\theta  )=0,
$$
and shows that the resulting estimator is not efficient for non-Markovian cases. In \cite{feuerverger_JCF} it is claimed that the variance of the estimator can be arbitrarily close to the Cramer-Rao bound if $n$ is chosen sufficiently large, but no proof is presented. Even if this claim were valid the implementation of the procedure for large $L$ would problematic. Moreover, Carrasco et all. argues that for large $n$ the available data provide only a few observation vectors of length $n.$

In this section as an alternative of the ML method we adapt the ECF method for GARCH processes. The motivation behind the adaptation of the ECF method again lies in the fact that the density function of $\Delta L_n$ is typically unknown. Still our proposed procedure estimates $\theta^*$ as efficiently as the ML method. We suppose that the characteristic $\eta^*$ of the noise is given and we are to identify the system parameters $\theta^*.$ Although the three-stage method can be applied for GARCH processes to identify the system and the noise characteristics, the results that we obtained for linear systems cannot be reproduced. The problem of identifying both the system parameters and the noise parameters will be briefly discussed at the end of the paper. The following paragraphs present the identification method with known $\eta^*.$

First, for each $\theta$ we define the estimated volatility $\sigma^2_n(\theta)$ and the estimated driving noise $\ino$ for $\theta \in D_{\epsilon},$ see equations (\ref{eq:sigma_inverse}) and (\ref{eq:vareps_def}). Following the philosophy of the ECF method take a fix set $u_i$-s, $1 \leq i \leq M.$ We define the $p \times 1$-dimensional modified primary score functions as
\begin{equation}\label{eq:garch_score}
h_{k,n}(\theta)=\left(e^{iu_k \varepsilon_n(\theta)}-\varphi(u)\right)\frac{\sigma_{\theta n}(\theta)}{\sigma_n(\theta)},
\end{equation}
where the modification being the usage of the instrumental variable $\frac{\sigma_{\theta n}(\theta)}{\sigma_n(\theta)}.$
The choice of the instrumental variable $\frac{\sigma_{\theta n}(\theta)}{\sigma_n(\theta)}$ is suggested by the construction of ECF method for linear systems. Namely, recall that for linear L\'evy systems the modified primary  score functions were defined via
$$
\left(e^{iu_k \varepsilon_n(\theta)}-\varphi(u)\right)\varepsilon_{\theta n}(\theta),
$$
 where the instrumental variable $\varepsilon_{\theta n}(\theta)$ satisfies $\lim_{n \rightarrow \infty} \E \left[\varepsilon_{\theta n}(\theta^*)\varepsilon^{T}_{\theta n}(\theta^*)\right]=R_P^*.$ By analogical thinking for GARCH processes we choose the instrumental variable $IV_n(\theta)$ such that $\lim_{n\rightarrow \infty}\E \left[IV_n(\theta^*)IV^{T}_n(\theta^*)\right]=M^*,$ hence the choice of $$\frac{\sigma_{\theta n}(\theta)}{\sigma_n(\theta)}.$$ Surprisingly we will see that this ad-hoc choice of instrumental variable yields an essentially asymptotically efficient identification method.

Since $\sigma_n(\theta),\sigma_{\theta n}(\theta)$ are $\mathscr{F}^{\Delta L}_{n-1}$ measurable
$$\E \left[h^{{\rm (s)}}_{k,n}(\theta^*)\right]=0$$
holds.
In analogy with the linear case merge the score functions $h_{k,n}(\theta)$-s into a $(r+s+1)M$-dimensional column vector
$$
h_n(\theta)=\left( h^T_{1,n}(\theta),\ldots,h^T_{M,n}(\theta)\right)^T.
$$
Define $\overline{h}_N(\theta)=\frac1N \sum_{n=1}^N h_{n}(\theta)$ the averaged score vector and
$$
g_N(\theta)=\E\left[ \overline{h}_N(\theta) \right] \text{~and~} g(\theta)=\lim_{N\rightarrow \infty} g_N(\theta).
$$
Note that the system of equations
$$
g(\theta)=0 %\quad n=1,\ldots,L
$$
is over-determined with solution $\theta=\theta^*,$ hence we redefine the score function as follows. Fix a symmetric, positive definite, $pM \times pM$ weighting matrix $K$. Since $g$ is not computable we approximate it by $\overline{h}_N$
and we seek a solution for the 'half-gradient' equation
\begin{equation}\label{eq:score2_2_garch}
V'_N(\theta)=\overline{h}^*_{\theta N}(\theta)K^{-1}\overline{h}_{N}(\theta)=0
\end{equation}
%uj{\bf !!! I WOULD NOT REPEAT EACH EQUATION BY PRESENTING ITS STATIONARY APPROXIMATION !}
to obtain $\hat{\theta}_N.$ We note in passing that the system of equations in (\ref{eq:score2_2_garch}) is no longer
over-determined because ${\rm dim }~V_{N}=r+s+1.$
Define
$$
G=g_{\theta}(\theta^*),
$$
and the auxiliary equation
\begin{equation}\label{eq:garch_auxiliary}
\bar{V}'_N(\theta)=G^*K^{-1}\overline{h}_{N}(\theta)=0.
\end{equation}
The asymptotic cost function is then given by
$$\bar{W}'(\theta)=\lim_{N \rightarrow \infty} \E[V'_{N}(\theta)]=g_{\theta}^*(\theta)K^{-1}g(\theta),$$
and its Jacobian at $\theta=\theta^*$ is
$$R_G^*=G^* K^{-1}G.$$
{\bf Condition}%\label{cond:garch_solutioncond:garch_solution}
The equation
$\bar{W}'(\theta)=0$ has a unique solution $\theta=\theta^*$ in
$D_{\epsilon}^*.$
%\end{condition}
We will use our recurring
$M \times M$ auxiliary matrix $C$ with elements
$$
C_{k,l} = \varphi(u_k - u_l,\eta^*) - \varphi(u_k ,\eta^*) \varphi(- u_l,\eta^*),
$$
recall that $C$ is the covariance matrix of the primary score functions used in the i.i.d. ECF method.

\section{Analysis of the ECF method for GARCH processes}

To analyze the process $\sigs$ we expand the state vector
$$
X^*_n=(y_n^2,\ldots,y^2_{n-r+1},\sigma_n^2,\ldots,\sigma_{n-s+1}^2)^T
$$ to
\begin{equation}\label{eq:state_garch_expanded}
\overline{X}_n(\theta)=(X_n^{*T},\sigs,\ldots,\sigma^2_{n-s+1}(\theta))^T.
\end{equation}
Then the dynamics of $\overline{X}_n(\theta)$ can be written as
\begin{equation}\label{eq:dynamics_garch_expanded}
\overline{X}_{n+1}(\theta)=\overline{A}_{n+1}(\theta) \overline{X}_n(\theta)+\overline{u}_{n+1}^*,
\end{equation}
where
$$
\overline{A}_{n}(\theta)=\left(
  \begin{array}{cc}
    A_n^* & Z \\
    M_{2,1}(\theta) & M_{2,2}(\theta) \\
  \end{array}
\right),
$$
with $Z$ being an $(r+s)\times s$ zero matrix,
$$
M_{2,1}(\theta)=\left(
                 \begin{array}{ccc|ccc}
                   \alpha_1 & \cdots & \alpha_r & 0 & \cdots & 0 \\
                   \hline
                   0 & \cdots & 0 & 0 & \cdots & 0 \\
                   \vdots &  \ddots & \vdots & \vdots & \ddots & \vdots \\
                   0 & \cdots & 0 & 0 & \cdots & 0 \\
                 \end{array}
               \right)
$$ is an $s \times (r+s)$ dimensional matrix, and
$$
M_{2,2}(\theta)=\left(
  \begin{array}{cccc|c}
    \beta_1 & \beta_2& \cdots & \beta_{s-1} & \beta_s \\
    \hline
    1 & 0 & \cdots & 0 & 0 \\
    0 & 1 & \cdots & 0 & 0 \\
    \vdots & \vdots & \ddots & \vdots & \vdots \\
    0 & 0 & \cdots & 1 & 0
  \end{array}
\right)
$$
is of dimension $s \times s,$ and finally
$$
\overline{u}_n^*=(u_n^{*T},0,\ldots,0)^T.
$$

First we state two theorems from the theory of block-triangular random matrices that we will use in the proofs, see \cite{orlovits_thesis}.
$\rho(P)$ stands for the spectral radius of matrix $P.$
\begin{theorem}\label{thm:garch_rho<1}
Let
$$
P=\left(
    \begin{array}{cc}
      P_1 & 0 \\
      B & P_2 \\
    \end{array}
  \right)
$$
be a random $(d_1+d_2) \times (d_1+d_2)$ matrix in $L^2(\Omega,\mathscr{F},P),$ with $P_1$ and $P_2$ being square matrices. Then
$$
\rho\left[ \E \left[ P \otimes P \right] \right]=\max \{ \rho\left[ \E \left[ P_1 \otimes P_1 \right] \right], \rho\left[ \E \left[ P_2 \otimes P_2 \right] \right]\}
$$
Similarly, let $q$ be a positive integer and let us assume that $P \in L^q(\Omega,\mathscr{F},P)$, then
$$
\rho\left[ \E \left[ P^{\otimes q} \right] \right]=\max \{ \rho\left[ \E \left[ P_1^{\otimes q} \right] \right]
,\rho\left[ \E \left[ P_2^{\otimes q} \right] \right]
\}.
$$
\end{theorem}

\begin{theorem}\label{thm:A_i product norm}
Let $(P_n)$ be an i.i.d. sequence of random matrices such that
$||P_1|| \in L^q.$ Assume that for some even integer $q \geq 2$
$$
\rho\left[ \E \left[ P_1^{\otimes q}  \right] \right] < 1
$$
holds. Then
$$
\lambda_q := \lim_{n \rightarrow \infty} \frac{1}{n} \log \E ||P_n \cdots P_1||^q <0.
$$
It follows that for any $\varepsilon > 0$ we have
$$
\E ||P_n \cdots P_1||^q \leq C e^{(\lambda_q+\varepsilon)n}
$$
with some $C=C(\varepsilon)>0$.
\end{theorem}
The next lemma implies the $L$-mixing property of the state vector.
\begin{lemma}
Let $D(q^{-1})$ be stable for all $\theta \in D_{\epsilon}$ and suppose that for some positive even $Q$ we have
$$
\rho\left[ \E\left[ (A_0^{*})^{\otimes Q} \right]\right]<1.
$$
Then the process $\left(\overline{X}_n(\theta)\right)$ is $L$-mixing of order Q uniformly in $\theta \in D_{\epsilon}$.
\end{lemma}

{\bf Proof:}
Fix a $\tau \in \mathbb{Z}^+$ and iterate the state space equation (\ref{eq:dynamics_garch_expanded})
\begin{equation}
\begin{split}
\overline{X}_n(\theta)=\overline{A}_n(\theta) \overline{X}_{n-1}(\theta) + \overline{u}_n^*=\overline{A}_n(\theta)\overline{A}_{n-1}(\theta)\overline{X}_{n-2}(\theta)+\overline{u}_n^*+ \overline{A}_{n}(\theta)\overline{u}_{n-1}^*= \ldots = \\
\overline{A}_n(\theta) \cdots \overline{A}_{n-\tau+1}(\theta) \overline{X}_{n-\tau}(\theta)+\overline{u}_n^*+\overline{A}_n(\theta)\overline{u}_{n-1}^*+ \ldots + \overline{A}_n(\theta) \cdots \overline{A}_{n-\tau+1}(\theta)\overline{u}_{n-\tau}^*
\end{split}
\end{equation}
Observe that
$$
\overline{u}_n^*+\overline{A}_n(\theta)\overline{X}_{n-1}(\theta)+ \ldots + \overline{A}_n(\theta) \cdots \overline{A}_{n-\tau+1}(\theta)\overline{u}_{n-\tau}^*
$$
is $\mathscr{F}_{n-\tau}^+=\sigma \{ \Delta L_i : i \geq n-\tau \}$ measurable, thus
\begin{equation*}
\begin{split}
&\E \left[ \overline{X}_n(\theta) | \mathscr{F}_{n-\tau}^+ \right]=
\overline{u}_n^*+\overline{A}_n(\theta)\overline{X}_{n-1}(\theta)+ \ldots + \overline{A}_n(\theta) \cdots \overline{A}_{n-\tau+1}(\theta)\overline{u}_{n-\tau}^*+ \\
&\E\left[ \overline{A}_n(\theta) \cdots \overline{A}_{n-\tau+1}(\theta) \overline{X}_{n-\tau}(\theta) | \mathscr{F}_{n-\tau}^+ \right]= \\
&\overline{u}_n^*+\overline{A}_n(\theta)\overline{X}_{n-1}(\theta)+ \ldots + \overline{A}_n(\theta) \cdots \overline{A}_{n-\tau+1}(\theta)\overline{u}_{n-\tau}^*+ \\
&\overline{A}_n(\theta) \cdots \overline{A}_{n-\tau+1}(\theta)\E\left[\overline{X}_{n-\tau}(\theta) \right],
\end{split}
\end{equation*}
because $\overline{X}_{n-\tau}(\theta)$ is independent of $\mathscr{F}_{n-\tau}^+.$
It follows that
\begin{equation}\label{eq:X_szetszedes}
\overline{X}_n(\theta)-\E \left[ \overline{X}_n(\theta) | \mathscr{F}_{n-\tau}^+ \right]=\overline{A}_n(\theta) \cdots \overline{A}_{n-\tau+1}(\theta) \left(\overline{X}_{n-\tau}(\theta)-\E\left[\overline{X}_{n-\tau}(\theta) \right]\right).
\end{equation}

Since $\overline{X}_{n-\tau}(\theta)$ is independent of $\overline{A}_n(\theta) \cdots \overline{A}_{n-\tau+1}(\theta)$ and $||AB||\leq ||A||~||B||$ for the $L^q$-norm of
(\ref{eq:X_szetszedes}) we have
\begin{equation}\begin{split}
\E^{1/q}\left[ ||\overline{A}_n(\theta) \cdots \overline{A}_{n-\tau+1}(\theta) \left(\overline{X}_{n-\tau}(\theta)-\E\left[\overline{X}_{n-\tau}(\theta) \right]\right)||^q\right] \leq  \\
\E^{1/q}\left[ ||\overline{A}_n(\theta) \cdots \overline{A}_{n-\tau+1}(\theta) ||^q\right] \E^{1/q}\left[|| \overline{X}_{n-\tau}(\theta)-\E\left[\overline{X}_{n-\tau}(\theta) \right]||^q\right].
\end{split}
\end{equation}
It is easy to see that $\overline{X}_{n}(\theta)-\E\left[\overline{X}_{n}(\theta) \right]$ is $M$-bounded of order $Q$, and for the first term of the two-term product on the l.h.s. using Theorem \ref{thm:garch_rho<1} with the choice $P_1=A_n^*$ and $P_2=M_{2,2}(\theta)$ yields $\rho\left[ \E \left[ \overline{A}_1(\theta)^{\otimes q}  \right] \right] < 1.$ Note that in this case the trivial version of Theorem \ref{thm:garch_rho<1} is used as $P_2$ is non-random. Hence, Theorem \ref{thm:A_i product norm} implies that
$$
\E^{1/q} ||\overline{A}_n(\theta) \cdots \overline{A}_{n-\tau+1}(\theta)||^q \leq C^{1/q} e^{(\lambda_q+\varepsilon)\tau/q}.
$$
Then choose $\varepsilon >0$ such that $\lambda_q+\varepsilon < 0.$ It follows that $\gamma_q (\tau,\overline{X}(\theta))$ is summable, which means by definition that $\left(\overline{X}(\theta)_n\right)$ is $L$-mixing or order $Q$ uniformly in $\theta \in D_{\epsilon}$. $\square$

\begin{lemma}
The process $\overline{X}_{e,n}(\theta):=\left(\overline{X}^T_n(\theta),\overline{X}^T_{\theta n}(\theta)\right)^T$ is $L$-mixing of order $Q$ uniformly in $\theta \in D_{\epsilon}$.
\end{lemma}
{\bf Proof:}
In order to analyze the derivative process we first determine its dynamics.
Suppose that we have a general parameter dependent recursion given by
\begin{equation}\label{eq:state_der_bevezeto}
\xi_{n+1}(\theta)=F_{n+1}(\theta) \xi_n(\theta)+v_{n+1}(\theta),
\end{equation}
and we are interested in the dynamic of the derivative process $\xi_{\theta n}(\theta).$ For simplicity we assume that $\theta$ is a scalar parameter, differentiating (\ref{eq:state_der_bevezeto}) we obtain
\begin{equation}
\xi_{\theta, n+1}(\theta)=F_{\theta, n+1}(\theta) \xi_n(\theta)+ F_{n+1}(\theta) \xi_{\theta,n}(\theta)+v_{\theta,n+1}(\theta).
\end{equation}
Thus the dynamics of the extended state vector $\xi_{e,n}=(\xi^T_n(\theta),\xi^T_{\theta,n}(\theta))^T$ can be written in a compact form:
\begin{equation}
\xi_{e,n+1}(\theta)=F_{e,n+1}(\theta) \xi_{e,n}(\theta)+v_{e,n+1}(\theta),
\end{equation}
with
$$
F_{e,n}(\theta)=\left(
                   \begin{array}{cc}
                     F_n (\theta)& 0 \\
                     F_{\theta,n} (\theta)& F_n (\theta)\\
                   \end{array}
                 \right),
$$
and $v_{e,n}(\theta)=\left(v^T_n(\theta),v^T_{\theta,n}(\theta)\right)^T.$% Using the following result we obtain $\rho\left[ \E \left[ \overline{A}_n(\theta)^{\otimes q} \right] \right].$

It follows that the state-transition matrix, say $\overline{A}_{e,n}(\theta),$ of the dynamics of $\overline{X}_{e,n}(\theta)$ has two identical blocks in the diagonal, namely $\overline{A}_n(\theta)$-s.
Hence Theorem \ref{thm:A_i product norm} implies that
$$
\rho\left[ \E \left[ \overline{A}_{e,n}(\theta)^{\otimes q} \right] \right]=\rho\left[ \E \left[ \overline{A}_n(\theta)^{\otimes q} \right] \right]<1.
$$
Mimicking the steps of the proof of the previous lemma we obtain that $\overline{X}_{e,n}(\theta)$ is $L$-mixing of order $Q$ uniformly in $\theta \in D_{\epsilon}$. Similarly, the same can be shown if we further expand $\overline{X}_{e,n}(\theta)$ with the higher order derivatives of $\sigma_{n}(\theta).$
$\square$

As consequence we get that
$$
h_{k,n}(\theta)=\left(e^{iu_k \varepsilon_n(\theta)}-\varphi(u)\right)\frac{\sigma_{\theta n}(\theta)}{\sigma_n(\theta)}
$$
and their derivatives w.r.t. $\theta$ up to order three are $L$-mixing of order $Q$. %From now on we may proceed as we did in Chapter\ref{chap:ECF_ARMA}.
%Define the averaged scores via
%$$
%\overline{h}_k(\theta)=\sum_{n=1}^{N} h_{k,n}(\theta),
%$$
%and merge these averaged errors into a vector $\overline{h}$
%$$
%\overline{h}(\theta)=\left(\overline{h}_k^T (\theta),\ldots,\overline{h}_M^T (\theta) \right)^T.
%$$
%Again introduce $g(\theta)$ to denote the expectation of $\overline{h}:$
%$$
%g(\theta)=\E \left[\overline{h} \right].
%$$
%We note here that $g$ is an $M$-dimensional vector.
We get the following major result which is a precise characterization of the estimation error:
\begin{theorem}
\label{th:difference_garch}
Assume that Condition 1 holds. Let $D(q^{-1})$ be stable for all $\theta \in D_{\epsilon}$ and suppose that for some positive even $Q$ we have
$$
\rho\left[ \E\left[ (A_0^{*})^{\otimes Q} \right]\right]<1.
$$
Then  for the estimation error we have
$$ \hat{\theta}_N-\theta^*=-(R_G^*)^{-1} \bar{V}'_{N}(\theta^*)+O^{Q/(2(r+s+1))}_M(N^{-1}).$$
\end{theorem}
The last formula equivalently can be written as
$$
\hat{\theta}_N-\theta^*=-(R_G^*)^{-1} G^*K^{-1}\overline{h}_N(\theta^*)+O^{Q/(2(r+s+1))}_M(N^{-1}).
$$
\section{Efficiency of the ECF method for GARCH processes}

In this section we show that the proposed ECF identification method gives an essentially asymptotically efficient estimate of the system characteristics of a GARCH process.
\begin{theorem}\label{thm:garch_asymp_cov}
Choose $K=C \otimes M^*,$ then for the estimate $\hat{\theta}_N$ obtained with the method presented in the previous section we have
$$
\E\left[N\left(\hat{\theta}_N-\theta^*\right)\left(\hat{\theta}_N-\theta^*\right)^*\right]=\Sigma_{\theta \theta}+O_M(N^{-1/2}),
$$
where the asymptotic covariance matrix is given by
$$\Sigma_{\theta \theta}=\left(\phi^* C^{-1} \phi\right)^{-1} (M^*)^{-1},$$
%uj{\bf !!! $\frac1N$ ???}
with $\phi=\left( u_1 \varphi' (u_1),\ldots,u_M \varphi'(u_M)\right)^T.$
\end{theorem}
The proof is analogous with that of Theorem 8 in \cite{automatica_own}. Note that $\phi$ and $\psi$ in the just mentioned theorem have similar structure, but now
\begin{equation*}
%\begin{split}
\E\left[h^{{\rm (s)}}_{\theta,k,n}(\theta^*)\right]=\E \left[e^{iu_k \Delta L_n} iu_k \varepsilon^{{\rm (s)}}_{\theta n}(\theta^*) \frac{\sigma^{{\rm (s)}}_{\theta n}(\theta^*)}{\sigma^{{\rm (s)}}_{n}(\theta^*)}\right]=
u_k\varphi'(u_k,\eta^*)M^*.
%\end{split}
\end{equation*}
Now we will demonstrate that the above presented ECF method gives an essentially asymptotically efficient estimate $\hat{\theta}_N.$ The line of reasoning is analogous with the one in the proof of Theorem 9 in \cite{automatica_own}. Suppose that we use the full continuum of moment conditions. Then the continuous version of (\ref{eq:garch_auxiliary}) would read as
$$
<K^{-1}G,\overline{h}_N>=0,
$$
where the inner product is defined on $H=L^2(\pi)=\left\{ f:\mathbb{R} \rightarrow \mathbb{C}\left|\int |f(t)|^2 \pi(t) dt < \infty \right.\right\}$ via
$$
<f,g>=\int f(t)g^*(t) \pi(t) dt,
$$
with $\pi$ being a probability measure on $\mathbb{R}.$

Define the $\pi$-dependent covariance operator %like in
%\cite{CARRASCO-EFFEMPIRCHAR}
%uj{\bf !!! WRITE OUT IT HERE !}
\begin{equation}
(Cf)(s)=\int c(s,t)f(t) \pi(t) dt,
\end{equation}
with
$$
c(s,t)=\E\left[
h_{s,n}(\theta^*,\eta^*)h^*_{t,n}(\theta^*,\eta^*)\right].
$$
If the full continuum of $u$-s were defined via $u_s=s$ for all
$s \in \mathbb{R},$ then the continuous version of Theorem \ref{thm:garch_asymp_cov} would give
\begin{equation}\label{eq:fullcont_cov}
\lim_{N \rightarrow \infty}\E\left[N\left(\hat{\theta}_N-\theta^*\right)\left(\hat{\theta}_N-\theta^*\right)^*\right]=\left(|| u\varphi'(u,\eta^*)||_C^2 \right)^{-1} (M^*)^{-1}
\end{equation}
for the asymptotic covariance matrix of the estimate $\hat{\theta}_N$.
Note that, like for linear system, in the above formula the asymptotic covariance matrix decouples, $|| u\varphi'(u,\eta^*)||_C^2$
depends only on $\eta^*$ and $R_P^*$ depends on the parameters
of the GARCH system.

Now we are ready to demonstrate that the proposed estimation
method is essentially asymptotically efficient provided the full continuum of moment conditions is available.
\begin{theorem}\label{thm:garch_efficiency}
Under the conditions of Theorem \ref{thm:garch_asymp_cov} the estimate $\hat{\theta}_N$ is essentially asymptotically efficient.
\end{theorem}
{\bf Proof:}
Recall that the asymptotic covariance of the ML estimate of the parameters of GARCH processes is
\begin{equation}
\mu^{-1} (M^*)^{-1},
\end{equation}
with
$$ M^*=\E\left[\frac{\sigma^{{\rm (s)}}_{\theta N}(\theta^*)\sigma^{{\rm (s)T}}_{\theta N}(\theta^*)}{\sigma^{{\rm (s)}2}_N(\theta^*)}\right],$$
and
\begin{equation}\label{eq:mu_def}
\mu=\E \left[ \frac{f'^2(\Delta L_n)}{f^2(\Delta L_n)}(\Delta L_n)^2-1\right].
\end{equation}
%Thus the ECF method gives an efficient estimate of the parameters of the GARCH process if
To complete the proof we only need to prove the following lemma.
\begin{lemma}
Using the notations above we have
\begin{equation}\label{eq:garch_mu_char}
\left(|| u\varphi'(u,\eta^*)||_C^2 \right)^{-1}=\mu^{-1}.
\end{equation}
\end{lemma}
Again, we do not prove (\ref{eq:garch_mu_char}) using direct computation. Instead
we show that $$\left(|| u\varphi'(u,\eta^*)||_C^2 \right)^{-1}$$
can be obtained as the asymptotic covariance of an efficient
ECF method with the full continuum moment conditions for the problem of estimating the scale parameter $\lambda^*$ of $\lambda^*\Delta L,$ with $\lambda^*=1,$ given an i.i.d. realization of $\Delta L.$  The problem of efficiency is then reduced to the i.i.d. case.

To carry out the suggested argument solve the following
identification problem: estimate the scale parameter $\lambda$ given a sequence of i.i.d. realizations of the
distribution $\lambda\Delta L$, where the true value of $\lambda$ is $\lambda=\lambda^*=1.$ The characteristic function of $\Delta Z$ is denoted by $\varphi,$ then the c.f. $\varphi_{\lambda \Delta L}(u,\lambda)$ of $\lambda \Delta L$ is given by $\varphi(u\lambda).$

Recall that for an i.i.d. sample, which was generated by a random variable
with a general characteristic function $\chi(u,\alpha^*),$ with $\alpha^*$ being an unknown parameter, the
ECF method using the full continuum of $u$-s gives an asymptotically efficient
estimate of $\alpha^*$ with asymptotic
covariance
$$(||\chi_{\alpha}(u,\alpha^*)||_C^2)^{-1}.$$

Write the derivative of the c.f. of $\lambda \Delta L$ w.r.t. $\lambda$
$$
\frac{\partial}{\partial \lambda}\E\left[e^{iu\lambda \Delta L}\right]= \\
\E\left[e^{iu \lambda \Delta L} iu\Delta L\right],
$$
choosing $\lambda=\lambda^*=1$ gives
$$
\left. \frac{\partial}{\partial \lambda}\varphi_{\lambda \Delta L}(u,\lambda) \right|_{\lambda=\lambda^*} =u \varphi'(u,\lambda^*) .
$$
Choosing $\chi=\varphi_{\lambda \Delta L}$ and $\alpha=\lambda$
we have $\chi_{\alpha}(u,\lambda^*)=u\varphi'(u,\eta^*).$
Hence for this identification problem the asymptotic covariance of the i.i.d. ECF method with full continuum $u$-s is
$$ \left(|| u\varphi'(u,\eta^*)||_C^2 \right)^{-1}.$$
Since the ECF method with continuum $u$-s is exactly as efficient as the ML method we find that
$ \left(|| u\varphi'(u,\eta^*)||_C^2 \right)^{-1}$
equals to the inverse Fisher of the ML method, hence (\ref{eq:garch_mu_char}) follows.
$\square$

\section{Discussion}

The optimal choice of $K$ is $C \otimes M^*,$ but $M^*$ is given by an expected value using the true value of parameters $\theta,$ so the optimal weighting matrix, like the optimal weighting matrix in the three-stage method for linear systems, is not computable.
We propose to approximate $M^*$ is two steps. First, define the approximation
$\hat{R}_P^*(\theta)$ by
$$
\hat{M}^*(\theta)=\frac1N\sum_{n=1}^{N} \frac{\sigma_{\theta n}(\theta)\sigma^{T}_{\theta n}(\theta)}{\sigma^2_n(\theta)}
$$
It would be convenient to use $\hat{M}^*(\theta^*),$ but since $\theta^*$ is unknown we approximate it by $\hat{\theta}^{{\rm (pre)}}_N,$ where $\hat{\theta}^{{\rm (pre)}}_N$ is a preliminary estimate obtained by using the ECF method for GARCH systems with the choice $K=I.$ Thus, we apply the ECF method with the weighting matrix
$$K=C \otimes \hat{M}^*\left(\hat{\theta}^{{\rm (pre)}}_N\right)$$ to get the approximation of $\hat{\theta}_N.$ It is relatively easy to see that Theorem \ref{th:difference_garch} and Theorem \ref{thm:garch_asymp_cov} are valid for this approximation of $\hat{\theta}_N,$ too.

As we have already mentioned at the beginning of the paper that although the three-stage method can be applied for GARCH processes to identify both the system and the noise characteristics, the results of \cite{automatica_own} cannot be reproduced. In what follows we address this issue.
Being aware of the steps of the three-stage method for linear L\'evy systems a three-stage identification method for GARCH systems can be proposed in a natural manner. Suppose now that both $\theta^*$ and $\eta^*$ are unknown.
The steps of the proposed three-method can be summarized as follows:
\begin{enumerate}
  \item Firstly estimate $\theta^*$ by applying the quasi-maximum likelihood method to obtain $\hat{\theta}_N.$
  \item Secondly invert the GARCH system with $\theta=\hat{\theta}_N$ to generate the estimated noise process,
  then estimate $\eta^*$ by pretending that these residuals are i.i.d., and apply the ECF method for i.i.d. data to obtain $\hat{\eta}_N.$
  \item Finally re-estimate $\theta^*$ by applying the ECF method for system identification, pretending that $\hat{\eta}_N = \eta^*,$
  to obtain an estimate $\hat{\hat{\theta}}_N$ for the dynamics.
\end{enumerate}
The problem with this three-stage method is that the tools presented for linear L\'evy systems cannot be adapted for its analysis. For, in analogy with the three-stage method for linear L\'evy systems the third step of the algorithm should give a consistent estimate of $\theta^*$ even if the noise characteristics $\eta$ is misspecified.
The $\eta$-dependent modified primary scores of the third step would be given by
\begin{equation*}
h_{k,n}(\theta,\eta)=\left(e^{iu_k \varepsilon_n(\theta)}-\varphi(u,\eta)\right)\frac{\sigma_{\theta n}(\theta)}{\sigma_n(\theta)}.
\end{equation*}
Following the notations and the line of arguments of Section \ref{sec:ECFmethod_GARCH} in defining the $\eta$-dependent scores, 'half-gradient' equations and corresponding variables
the asymptotic value of function $V'_{N}(\theta,\eta)$ would be given by
$$\bar{W}'(\theta,\eta)=\lim_{N \rightarrow \infty} \E[V'_{N}(\theta,\eta)]=g_{\theta}^*(\theta,\eta)K^{-1}g(\theta,\eta).$$
Observe that if we are given a misspecified $\eta,$ then by solving $\bar{W}'(\theta,\eta)=0$ for $\theta$ we typically have a solution $\theta^*(\eta)$ such that
$$
\theta^*(\eta) \neq \theta^*.
$$
The reason behind is that for the instrumental variable we typically have that
$$
\E\left[ \frac{\sigma_{\theta n}(\theta^*)}{\sigma_n(\theta^*)}\right] \neq 0.
$$
The study of this interesting problem will be a subject of our further research.

\bibliography{refereces_new}{}
\bibliographystyle{abbrv}
%\appendix
%\section{Proofs}    % Each appendix must have a short title.
%\section{Some Latin vocabulary}         % Sections and subsections are supported
                                        % in the appendices.
\end{document}